\documentclass[preprint,12pt]{elsarticle}

\usepackage{amssymb}
 \usepackage{amsthm}
\usepackage{latexsym}
\usepackage{amsmath}
\usepackage{amssymb}
\usepackage{slashbox}
\usepackage{arydshln}

\newcommand{\udots}{\mathinner{\mskip1mu\raise1pt\vbox{\kern7pt\hbox{.}}
\mskip2mu\raise4pt\hbox{.}\mskip2mu\raise7pt\hbox{.}\mskip1mu}}
\usepackage{amsmath}

\newtheorem{theorem}{Theorem}

\newtheorem{lemma}{Lemma}
\newtheorem{remark}{Remark}
\newtheorem{proposition}{Proposition}

\begin{document}

\begin{frontmatter}



\title{Expressions for joint moments of elliptical distributions
}


\author[label1]{Baishuai Zuo}
\author[label1]{Chuancun Yin\corref{cor1}}
 \ead{ccyin@qfnu.edu.cn}
 \author[label2]{Narayanaswamy Balakrishnan\corref{cor1}}
 \cortext[cor1]{Corresponding author.}
 \ead{bala@mcmaster.ca}
\address[label1]{ School of Statistics, Qufu Normal University, Qufu, Shandong 273165, China}
\address[label2]{ Department of Mathematics and Statistics, McMaster University, Hamilton, Ontario, Canada
}
\begin{abstract}
 Inspired by Stein's lemma,  we derive two expressions for the  joint moments of elliptical distributions. We use two different methods  to derive $E[X_{1}^{2}f(\mathbf{X})]$ for any measurable function $f$ satisfying some regularity conditions. Then, by applying this result, we obtain new formulae for expectations of product of normally distributed random variables, and also present simplified expressions of $E[X_{1}^{2}f(\mathbf{X})]$ for multivariate Student-$t$, logistic and Laplace distributions.
\end{abstract}

\begin{keyword}

 Multivariate elliptical distributions;   Multivariate normal distribution; Joint moments; Stein's lemma
\end{keyword}

\end{frontmatter}

\baselineskip=20pt

\section{Introduction and motivation}
Stein's (1981) lemma discusses the determination of $Cov(X, h(Y))$  for a bivariate normal random vector $(X,Y)$,
 where $h$ is any differentiable function with finite extensions $E|h'(Y)|$.  Inspired by  the original work of Stein, several and generalizations
have appeared in the literature. Liu (1994) generalized the lemma to the multivariate normal case, while
Landsman (2006) showed that Stein's type lemma also holds when $(X,Y)$ is distributed as bivariate elliptical. This result got extended to multivariate elliptical vectors by Landsman
and Neslehova (2008); see also Landsman et al.(2013) for a simple proof. Recently, Shushi (2018) derived the multivariate Stein's lemma for truncated elliptical random vectors.

 In this work,  we  derive expressions of the joint moments $E[X_{1}^{2}f(\mathbf{X})]$ for any measurable function $f$ satisfying some regularity conditions.  In particular,  we obtain new formulae for the expectation of product of normally distributed random variables, and also present expressions of $E[X_{1}^{2}f(\mathbf{X})]$ for multivariate Student-$t$, logistic and Laplace distributions.

The rest of the paper is organized as follows. Section 2 reviews some definitions and properties of the family of elliptical distributions. Section 3 presents explicit expressions of joint moments of elliptical distributions by a direct method. Section 4 derives these joint moments by the use of Stein's lemma. Section 5 shows the equivalence of the two expressions under the condition that the  scale matrix  is positive definite. Section 6 presents expression for expectation of products of correlated normal variables. Section 7 presents simplified results for the special cases of multivariate Student-$t$, logistic and Laplace distributions as illustrative examples of the general result established here. Finally, Section 8 gives some concluding remarks.
\section{Family of elliptical distributions}

Elliptical distributions are generalizations of multivariate normal distribution and possess many novel tractable properties, in addition to allowing fat tails with suitably chosen keernels. This class of distributions was first introduced by Kelker (1970), and has been widely discussed in detail by Fang et al. (1990), and Kotz et al. (2000). A $n \times 1$ random vector $X = (X_1,\cdots, X_n)^{T}$ is said to have an elliptically
symmetric distribution if its characteristic function has the form
 $$E[\exp(i{\bf t}^{T}{\bf X})]= e^{i{\bf t}^{T}{\boldsymbol \mu}}\phi\left(\frac12 {\bf t}^{T}{\bf \Sigma}{\bf t}\right)$$
  for all ${\bf t}\in \Bbb{R}^n $,  denoted ${\bf{X}}\sim E_n ({\boldsymbol \mu},{\bf \Sigma},\phi)$,
  where  $\phi$ is called the characteristic generator with $\phi(0)=1$,
$\boldsymbol{\mu}$ ($n$-dimensional vector) is the location parameter,  and  $\bf{\Sigma}$ ($n\times n$ matrix with $\bf{\Sigma}\ge 0$) is the dispersion matrix (or scale matrix). The mean vector $E({\bf X})$   (if it exists)
coincides with the location vector and the covariance matrix  Cov$({\bf X})$ (if it exists) is $-\phi'(0){\bf \Sigma}$.
  The generator of the multivariate normal distribution, for example, is given by $\phi(u)=\exp(-u)$.

In general, the elliptical  vector ${\bf{X}}\sim E_n ({\boldsymbol \mu},{\bf \Sigma},\phi)$ may not have a density function. However, if the density
 $f_{\bf X}(\boldsymbol{x})$  exists, then it is of the form
\begin{align}\label{(P1)}
f_{\boldsymbol{X}}(\boldsymbol{x})=\frac{c_{n}}{\sqrt{|\mathbf{\Sigma}|}}
g_{n}\left(\frac{1}{2}(\boldsymbol{x}-\boldsymbol{\mu})^{T}\mathbf{\Sigma}^{-1}(\boldsymbol{x}-\boldsymbol{\mu})\right),\;  \boldsymbol{x}\in \Bbb{R}^n,
\end{align}
where $\boldsymbol{\mu}$ is an $n\times1$ location vector, $\mathbf{\Sigma}$ is an $n\times n$ positive definite scale matrix, and $g_{n}(u)$, $u\geq0$, is the density generator of $\mathbf{X}$. This density generator satisfies the condition
\begin{align*}
\int_{0}^{\infty}t^{(n/2)-1}g_{n}(t)\mathrm{d}t<\infty,
\end{align*}
 and the normalizing constant  $c_n$ is given by
\begin{align}\label{(P2)}
c_{n}=\frac{\Gamma(n/2)}{(2\pi)^{n/2}}\left[\int_{0}^{\infty}t^{(n/2)-1}g_{n}(t)\mathrm{d}t\right]^{-1}.
\end{align}
Two important special cases are multivariate normal family with $g_n(u)=e^{-u}$, and
multivariate generalized Student-$t$ family with $g_n(u)=(1+\frac{u}{k_{n,p}})^{-p}$, where the parameter $p>\frac{n}{2}$ and $k_{n,p}$ is some
constant that may depend on $n$ and $p$.

To derive the mixed moments of elliptical distributions, we use the cumulative generators $\overline{G}_{n}(u)$ and $\overline{\mathcal{G}}_{n}(u)$, which are given by
\begin{align}\label{(P3)}
\overline{G}_{n}(u)=\int_{u}^{\infty}{g}_{n}(v)\mathrm{d}v
\end{align}
and
\begin{align}\label{(P4)}
\overline{\mathcal{G}}_{n}(u)=\int_{u}^{\infty}{\overline{G}}_{n}(v)\mathrm{d}v,
\end{align}
respectively (see Landsman et al. (2018)), and the corresponding normalizing constants  are
\begin{align}\label{(P5)}
c_{n}^{\ast}=\frac{\Gamma(n/2)}{(2\pi)^{n/2}}\left[\int_{0}^{\infty}t^{n/2-1}\overline{G}_{n}(t)\mathrm{d}t\right]^{-1}
\end{align}
and
\begin{align}\label{(P6)}
c_{n}^{\ast\ast}=\frac{\Gamma(n/2)}{(2\pi)^{n/2}}\left[\int_{0}^{\infty}t^{n/2-1}\overline{\mathcal{G}}_{n}(t)\mathrm{d}t\right]^{-1}.
\end{align}

Throughout this paper, ${\bf x}\in\Bbb{R}^n$ will denote an $n$-dimensional vector and ${\bf x}^{T} =(x_1,\cdots, x_n )$ its transpose.
For an $n \times n$ matrix ${\bf \Sigma}\in \Bbb{R}^{n\times n}$, $|{\bf \Sigma}|$ is the determinant of ${\bf \Sigma}$. If ${\bf \Sigma}$ is positive definite, then its Cholesky decomposition is known to be unique.

\section{Direct method of derivation}
  Consider a random vector $\mathbf{X}\sim E_{n}(\boldsymbol{\mu},~\mathbf{\Sigma},~g_{n})$ with mean vector $\boldsymbol{\mu}=(\mu_1,\cdots,\mu_n)^{T}$ and positive define matrix ${\bf {\Sigma}}= (\sigma_{ij})_{i,j=1}^{n}$.
   Partition $\boldsymbol{y}\in \Bbb{R}^n$ into two parts as $\boldsymbol{y}=(y_{1},\boldsymbol{y}_{(2)})^{T}$ each with $1$ and $n-1$ components, respectively.
 By Cholesky decomposition (see Golub and Van Loan (2012), for example),   there  exists a unique lower triangular matrix $\mathbf{A}=(a_{ij})_{i,j=1}^{n} $ such that $\mathbf{AA}^{T}=\mathbf{\Sigma}$. In terms of components, we have
  \begin{align}\label{(P7)}
&a_{11}=\sqrt{\sigma_{11}},~a_{i1}=\frac{\sigma_{i1}}{a_{11}},~i=2,\cdots,~n, ~a_{ik}=0, i<k,\\
&a_{kk}=\sqrt{\sigma_{kk}-\sum_{i=1}^{k-1}a^2_{ki}},~a_{ik}=\frac{\sigma_{ik}-\sum_{j=1}^{k-1}a_{ij}a_{kj}}{a_{kk}},~i=k+1,~\cdots,~n.\label{(P8)}
\end{align}

Let $f:\mathbb{R}^{n}\rightarrow \mathbb{R}$  be a twice continuously differentiable function, and shall use
   $(\nabla_{i,j}f( \boldsymbol{x}))_{i,j=1}^{n}$ to denote the Hessian matrix of $f$. In addition, we denote
\begin{align*}
\nabla_{i,j}f( \boldsymbol{x})=\frac{\partial^{2} f( \boldsymbol{x})}{\partial x_{i}\partial x_{j}},~i,j=1,~2,\cdots,n,
\end{align*}
\begin{align*}
\nabla_{i}f( \boldsymbol{x})=\frac{\partial f( \boldsymbol{x})}{\partial y_{i}},~i=1,~2,\cdots,n
\end{align*}
and
\begin{align*}
\nabla f(\boldsymbol{x})=\left(\frac{\partial f(\boldsymbol{x})}{\partial x_{1}},\frac{\partial f(\boldsymbol{x})}{\partial x_{2}},\cdots,\frac{\partial f(\boldsymbol{x})}{\partial x_{n}}\right)^{T}.
\end{align*}
Let $\mathbf{X}^{\ast}\sim E_{n}(\boldsymbol{\mu},~\mathbf{\Sigma},~\overline{G}_{n})$ and  $\mathbf{X}^{\ast\ast}\sim E_{n}(\boldsymbol{\mu},~\mathbf{\Sigma},~\overline{\mathcal{G}}_{n})$ be  two elliptical random vectors with generators $\overline{G}_{n}(u)$ and $\overline{\mathcal{G}}_{n}(u)$, repectively.

The following theorem  gives an  expression for joint moments of elliptical distributions.
\begin{theorem}\label{th.1}
Let $\mathbf{X}\sim E_{n}(\boldsymbol{\mu},~\mathbf{\Sigma},~g_{n})$ be an $n$-dimensional elliptical random vector with density generator $g_{n}$, positive definite matrix $\mathbf{\Sigma}=(\sigma_{i,j})_{i,j=1}^{n}$, and finite expectation $\boldsymbol{\mu}$. Further, let $f:\mathbb{R}^{n}\rightarrow \mathbb{R}$  be a twice continuously differentiable function satisfying $E[\nabla_{i,j}f(\boldsymbol{X}^{\ast\ast})]<\infty$ and $E[\nabla_{i}f(\boldsymbol{X}^{\ast})]<\infty$. Let, in addition,
\begin{align}\label{(P9)}
\lim_{|x_{1}|\rightarrow\infty}x_{1}f(\mathbf{A}\boldsymbol{x}+\boldsymbol{\mu})\overline{G}_{n}\left(\frac{1}{2}\boldsymbol{x}^{T}\boldsymbol{x}\right)=0
\end{align}
and
\begin{align}\label{(P10)}
\lim _{|x_{1}|\rightarrow\infty}[\nabla_{1} f(\mathbf{A}\boldsymbol{x}+\boldsymbol{\mu})]\overline{\mathcal{G}}_{n}\left(\frac{1}{2}\boldsymbol{x}^{T}\boldsymbol{x}\right)=0.
\end{align}
 Then,
\begin{align}\label{(P11)}
\nonumber E[X_{1}^{2}f(\mathbf{X})]&= \sigma_{11}b_{n}^{\ast}E[f(\mathbf{X}^{\ast})]+
b_{n}^{\ast\ast}\sum_{i=1}^{n}\sum_{j=1}^{n}\sigma_{i1}\sigma_{j1}E[\nabla_{i,j}f(\mathbf{X}^{\ast\ast})]\\
&+2\mu_{1}b_{n}^{\ast}\sum_{i=1}^{n}\sigma_{i1}E[\nabla_{i} f(\mathbf{X}^{\ast})]+\mu_{1}^{2}E[f(\mathbf{X})],
\end{align}
where   $b_{n}^{\ast}=\frac{c_{n}}{c_{n}^{\ast}}$ and $b_{n}^{\ast\ast}=\frac{c_{n}}{c_{n}^{\ast\ast}}$.
\end{theorem}
\noindent Proof. By definition,
\begin{align*}
E[X_{1}^{2}f(\mathbf{X})]&=\frac{c_{n}}{\sqrt{|\mathbf{\Sigma}|}}\int_{\mathbb{R}^{n}}x_{1}^{2}f(\boldsymbol{x})
g_{n}\left\{\frac{1}{2}(\boldsymbol{x}-\boldsymbol{\mu})^{T}\Sigma^{-1}(\boldsymbol{x}-\boldsymbol{\mu})\right\}\mathrm{d}\boldsymbol{x}\\
&=\frac{c_{n}}{\sqrt{|\mathbf{\Sigma}|}}\int_{\mathbb{R}^{n}}x_{1}^{2}f(\boldsymbol{x})
g_{n}\left\{\frac{1}{2}(\boldsymbol{x}-\boldsymbol{\mu})^{T}(\mathbf{AA}^{T})^{-1}(\boldsymbol{x}-\boldsymbol{\mu})\right\}\mathrm{d}\boldsymbol{x}.
\end{align*}
Now, setting $\boldsymbol{y}=\mathbf{A}^{-1}(\boldsymbol{x}-\boldsymbol{\mu})$, we obtain
\begin{align*}
E[X_{1}^{2}f(\mathbf{X})]=&\frac{c_{n}|\mathbf{A}|}{\sqrt{|\mathbf{\Sigma}|}}
\int_{\mathbb{R}^{n}}(a_{1,1}y_{1}+\mu_{1})^{2}f(\mathbf{A}\boldsymbol{y}+\boldsymbol{\mu})g_{n}\left\{\frac{1}{2}\boldsymbol{y}^{T}\boldsymbol{y}\right\}\mathrm{d}\boldsymbol{y}\\
=&c_{n}\int_{\mathbb{R}^{n}}(a_{1,1}y_{1}+\mu_{1})^{2}f(\mathbf{A}\boldsymbol{y}+\boldsymbol{\mu})
g_{n}\left\{\frac{1}{2}\boldsymbol{y}^{T}\boldsymbol{y}\right\}\mathrm{d}\boldsymbol{y}\\
=&a_{1,1}^{2}c_{n}\int_{\mathbb{R}^{n}}y_{1}^{2}f(\mathbf{A}\boldsymbol{y}+\boldsymbol{\mu})
g_{n}\left\{\frac{1}{2}\boldsymbol{y}^{T}\boldsymbol{y}\right\}\mathrm{d}\boldsymbol{y}\\
&+2a_{1,1}\mu_{1}c_{n}\int_{\mathbb{R}^{n}}y_{1}f(\mathbf{A}\boldsymbol{y}+\boldsymbol{\mu})
g_{n}\left\{\frac{1}{2}\boldsymbol{y}^{T}\boldsymbol{y}\right\}\mathrm{d}\boldsymbol{y}\\
&+\mu_{1}^{2}c_{n}\int_{\mathbb{R}^{n}}f(\mathbf{A}\boldsymbol{y}+\boldsymbol{\mu})
g_{n}\left\{\frac{1}{2}\boldsymbol{y}^{T}\boldsymbol{y}\right\}\mathrm{d}\boldsymbol{y}\\
=&a_{1,1}^{2}c_{n}\int_{\mathbb{R}^{n-1}}I_{1}d\boldsymbol{y}_{(2)}+2a_{1,1}\mu_{1}c_{n}\int_{\mathbb{R}^{n-1}}I_{2}\mathrm{d}\boldsymbol{y}_{(2)}\\
&+\mu_{1}^{2}c_{n}\int_{\mathbb{R}^{n}}f(\mathbf{A}\boldsymbol{y}+\boldsymbol{\mu})
g_{n}\left\{\frac{1}{2}\boldsymbol{y}^{T}\boldsymbol{y}\right\}\mathrm{d}\boldsymbol{y},
\end{align*}
where
\begin{align*}
I_{1}&=\int_{\mathbb{R}}y_{1}^{2}f(\mathbf{A}\boldsymbol{y}+\boldsymbol{\mu})g_{n}\left\{\frac{1}{2}\boldsymbol{y}^{T}\boldsymbol{y}\right\}\mathrm{d}y_{1}\\
&=-\int_{\mathbb{R}}y_{1}f(\mathbf{A}\boldsymbol{y}+\boldsymbol{\mu})\frac{\partial }{\partial y_1}\overline{G}_{n}\left\{\frac{1}{2}\boldsymbol{y}^{T}\boldsymbol{y}\right\}\mathrm{d}y_{1} \\
&=\int_{\mathbb{R}}[f(\mathbf{A}\boldsymbol{y}+\boldsymbol{\mu})+y_{1}\nabla_{1}f(\mathbf{A}\boldsymbol{y}+\boldsymbol{\mu})]
\overline{G}_{n}\left\{\frac{1}{2}\boldsymbol{y}^{T}\boldsymbol{y}\right\}\mathrm{d}y_{1}\\
&=\int_{\mathbb{R}}f(\mathbf{A}\boldsymbol{y}+\boldsymbol{\mu})
\overline{G}_{n}\left\{\frac{1}{2}\boldsymbol{y}^{T}\boldsymbol{y}\right\}\mathrm{d}y_{1}-\int_{\mathbb{R}}\nabla_{1}f(\mathbf{A}\boldsymbol{y}
+\boldsymbol{\mu})\frac{\partial}{\partial y_1}\overline{\mathcal{G}}_{n}\{\frac{1}{2}\boldsymbol{y}^{T}\boldsymbol{y}\}\mathrm{d}y_{1}\\
&=\int_{\mathbb{R}}f(\mathbf{A}\boldsymbol{y}+\boldsymbol{\mu})
\overline{G}_{n}\left\{\frac{1}{2}\boldsymbol{y}^{T}\boldsymbol{y}\right\}\mathrm{d}y_{1}+\int_{\mathbb{R}}\nabla_{1,1}f(\mathbf{A}\boldsymbol{y}
+\boldsymbol{\mu})\overline{\mathcal{G}}_{n}\left\{\frac{1}{2}\boldsymbol{y}^{T}\boldsymbol{y}\right\}\mathrm{d}y_{1},
\end{align*}
 and
\begin{align*}
I_{2}&=\int_{\mathbb{R}}y_{1}f(\mathbf{A}\boldsymbol{y}+\boldsymbol{\mu})g_{n}\left\{\frac{1}{2}\boldsymbol{y}^{T}\boldsymbol{y}\right\}\mathrm{d}y_{1}\\
&=-\int_{\mathbb{R}}f(\mathbf{A}\boldsymbol{y}+\boldsymbol{\mu}) \frac{\partial}{\partial y_1}\overline{G}_{n}\left\{\frac{1}{2}\boldsymbol{y}^{T}\boldsymbol{y}\right\}\mathrm{d}y_{1}\\
&=\int_{\mathbb{R}}\nabla_{1}f(\mathbf{A}\boldsymbol{y}+\boldsymbol{\mu})\overline{G}_{n}\left\{\frac{1}{2}\boldsymbol{y}^{T}\boldsymbol{y}\right\}
\mathrm{d}y_{1}.
\end{align*}
 We then obtain
\begin{align*}
E[X_{1}^{2}f(\mathbf{X})]&=a_{11}^{2}\{b_{n}^{\ast}E[f(\mathbf{X}^{\ast})]+b_{n}^{\ast\ast}\mathbf{A}_{1}^{T}E(\nabla_{i,j}f(\mathbf{X}^{\ast\ast}))_{i,j=1}^{n}\mathbf{A}_{1}\}\\
&+2a_{11}\mu_{1}b_{n}^{\ast}\mathbf{A}_{1}^{T}E[\nabla f(\mathbf{X}^{\ast})]+\mu_{1}^{2}E[f(\mathbf{X})],
\end{align*}
where $\mathbf{A}=(a_{ij})_{i,j=1}^{n}=(\mathbf{A}_{1},\cdots,\mathbf{A}_{n})$. Now, upon using (\ref{(P7)}) and (\ref{(P8)}),
we obtain (\ref{(P11)}), completing the proof of the theorem.

 Kan (2008) and Song and Lee (2015) have presented  explicit formulae for product moments of multivariate Gaussian random variables.
 The following corollary presents an explicit expression for product moments of multivariate elliptical random variables, in general.

 \noindent $\mathbf{Corollary~1}$. Suppose ${\bf{X}}\sim E_n ({\boldsymbol \mu},{\bf \Sigma},g_{n})$,  ${\bf{X}^*}\sim E_n ({\boldsymbol \mu},{\bf \Sigma},\overline{G}_{n})$ and  ${\bf{X}^{**}}\sim E_n ({\boldsymbol \mu},{\bf \Sigma},\overline{\mathcal{G}}_{n})$. Let
  $p_{1}, \cdots, p_{n}$ be nonnegative~integers with $p_{1}\geq2$.  Then, we have
\begin{align*}
\nonumber &E\left[\prod_{i=1}^{n}X_{i}^{p_{i}}\right]=b_{n}^{\ast}\sigma_{11}E\left[{(X^*_{1})}^{p_{1}-2}\prod_{k=2}^{n}{(X^*_{k})}^{p_{k}}\right]\\
\nonumber &+b_{n}^{\ast\ast}\bigg\{\sigma_{11}^{2}(p_{1}-2)(p_{1}-3)E\left[{(X^{**}_{1})}^{p_{1}-4}\prod_{k=2}^{n}{(X^{**}_{k})}^{p_{k}}\right]\\
\nonumber &+2\sum_{j=2}^{n}\sigma_{11}\sigma_{j1}(p_{1}-2)p_{j}E\left[{(X^{**}_{1})}^{p_{1}-3}{(X^{**}_{j})}^{p_{j}-1}\prod_{k=2,~k\neq j}^{n}{(X^{**}_{k})}^{p_{k}}\right]\\
\nonumber &+\sum_{j=2}^{n}\sigma_{j1}^{2}p_{j}(p_{j}-1)E\left[{(X^{**}_{1})}^{p_{1}-2}{(X^{**}_{j})}^{p_{j}-2}\prod_{k=2,~k\neq j}^{n}{(X^{**}_{k})}^{p_{k}}\right]\\
\end{align*}
\begin{align}\label{F17}
\nonumber &+\sum_{j=2}^{n}\sum_{i=2,~i\neq j}^{n}\sigma_{j1} \sigma_{i1}p_{j}p_{i}E\left[{(X^{**}_{i})}^{p_{i}-1}{(X^{**}_{j})}^{p_{j}-1}{(X^{**}_{1})}^{p_{1}-2}\prod_{k=2,~k\neq i,~j}^{n}{(X^{**}_{k})}^{p_{k}}\right]\bigg\}\\
\nonumber &+2b_{n}^{\ast}\mu_{1}\bigg\{\sigma_{11}(p_{1}-2)E\left[{(X^{**}_{1})}^{p_{1}-3}\prod_{k=2}^{n}{(X^{**}_{k})}^{p_{k}}\right]\\
&+\sum_{j=2}^{n}\sigma_{j1}p_{j}E\left[{(X^*_{j})}^{p_{j}-1}{(X^*_{1})}^{p_{1}-2}\prod_{k=2,~k\neq j}^{n}{(X^*_{k})}^{p_{k}}\right]\bigg\}+\mu_{1}^{2}E\left[X_{1}^{p_{1}-2}\prod_{k=2}^{n}X_{k}^{p_{k}}\right].
\end{align}

From Corollary~1, we readily deduce the following relations for the special case of $p_{1}=3$, for example:
\begin{align*}
\nonumber &E\left[X_{1}^{3}\prod_{i=2}^{n}X_{i}^{p_{i}}\right]=b_{n}^{\ast}\sigma_{11}E\left[X^*_{1}\prod_{k=2}^{n}{(X^*_{k})}^{p_{k}}\right]+\mu_{1}^{2}E\left[X_{1}\prod_{k=2}^{n}X_{k}^{p_{k}}\right]\\
&+b_{n}^{\ast\ast}\bigg\{2\sum_{j=2}^{n}\sigma_{11}\sigma_{j1}p_{j}E\left[{(X^{**}_{j})}^{p_{j}-1}\prod_{k=2,~k\neq j}^{n}{(X^{**}_{k})}^{p_{k}}\right]\\
&+\sum_{j=2}^{n}\sigma_{j1}^{2}p_{j}(p_{j}-1)E\left[X^{**}_{1}{(X^{**}_{j})}^{p_{j}-2}\prod_{k=2,~k\neq j}^{n}{(X^{**}_{k})}^{p_{k}}\right]\\
&+\sum_{j=2}^{n}\sum_{i=2,~i\neq j}^{n}\sigma_{j1} \sigma_{i1}p_{j}p_{i}E\left[{(X^{**}_{i})}^{p_{i}-1}{(X^{**}_{j})}^{p_{j}-1}X^{**}_{1}\prod_{k=2,~k\neq i,~j}^{n}{(X^{**}_{k})}^{p_{k}}\right]\bigg\}\\
&+2b_{n}^{\ast}\mu_{1}\bigg\{\sigma_{11}E\left[\prod_{k=2}^{n}{(X^{**}_{k})}^{p_{k}}\right]+\sum_{j=2}^{n}\sigma_{j1}p_{j}E\left[{(X^*_{j})}^{p_{j}-1}X^*_{1}\prod_{k=2,~k\neq j}^{n}{(X^*_{k})}^{p_{k}}\right]\bigg\}.
\end{align*}
\section{Derivation by the use of Stein's lemma}
We now derive an expression for $E[X_{1}^{2}f(\boldsymbol{X})]$ by using Stein's Lemma. It should be noted that
the  positive definiteness of ${\bf \Sigma}$ is not necessary in this case.
\begin{theorem}\label{th.2} Suppose ${\bf{X}}\sim E_n ({\boldsymbol \mu},{\bf \Sigma},\phi)$, and all the conditions of Theorem \ref{th.1} hold. Then,\\
\begin{align}\label{(PP1)}
\nonumber E[X_{1}^{2}f(\boldsymbol{X})]=&\sum_{i=1}^{n}\sum_{j=1}^{n}Cov(X_{1},X_{i})Cov(X_{1}^{\ast},X_{j}^{\ast})E[\nabla_{i,j}f(\boldsymbol{X}^{\ast\ast})]\\
&+2\mu_1\sum_{i=1}^{n}Cov(X_{1},X_{i})E[\nabla_{i}f(\mathbf{X}^{\ast})]\nonumber\\
&+Cov(X_{1},X_{1})E[f({\bf X}^{\ast})]+\mu_{1}^{2}E[f(\boldsymbol{X})].
\end{align}
\end{theorem}
\noindent Proof. By Lemma 2 of Landsman et al. (2013), we have
\begin{align*}
Cov(X_{1},f(\mathbf{X}))=\sum_{i=1}^{n}Cov(X_{1},X_{i})E[\nabla_{i}f(\mathbf{X}^{\ast})],
\end{align*}
and so
\begin{align}\label{(PP2)}
\nonumber E[X_{1}f(\boldsymbol{X})]&=Cov(X_{1},f(\mathbf{X}))+E[X_{1}]E[f(\mathbf{X})]\\
&=\sum_{i=1}^{n}Cov(X_{1},X_{i})E[\nabla_{i}f(\mathbf{X}^{\ast})]+E[X_{1}]E[f(\mathbf{X})].
\end{align}
Upon replacing $f(\mathbf{X})$ by  $X_{1}f(\mathbf{X})$  in  $(\ref{(PP2)})$, we obtain
 \begin{align*}
\nonumber &E[X_{1}^{2}f(\boldsymbol{X})]=\sum_{i=1}^{n}Cov(X_{1},X_{i})E[\nabla_{i}(X_{1}^{\ast}f(\mathbf{X}^{\ast}))]+E[X_{1}]E[X_{1}f(\mathbf{X})]\\
=&\sum_{i=1}^{n}Cov(X_{1},X_{i})E[X_{1}^{\ast}\nabla_{i}f(\mathbf{X}^{\ast})]+Cov(X_{1},X_{1})E[f(\mathbf{X}^{\ast})]+E[X_{1}]E[X_{1}f(\mathbf{X})]\\
=&\sum_{i=1}^{n}Cov(X_{1},X_{i})\left\{\sum_{j=1}^{n}Cov(X_{1}^{\ast},X_{j}^{\ast})E[\nabla_{i,j}f(\mathbf{X}^{\ast\ast})]+E[X_{1}^{\ast}]E[\nabla_{i}f(\mathbf{X}^{\ast})]\right\}\\
&+Cov(X_{1},X_{1})E[f(\mathbf{X}^{\ast})]+E[X_{1}]\left\{\sum_{i=1}^{n}Cov(X_{1},X_{i})E[\nabla_{i}f(\mathbf{X}^{\ast})]+E[X_{1}]E[f(\mathbf{X})]\right\}
\end{align*}
\begin{align*}
=&\sum_{i=1}^{n} \sum_{j=1}^{n} Cov(X_{1},X_{i})Cov(X_{1}^{\ast},X_{j}^{\ast})E[\nabla_{i,j}f(\mathbf{X}^{\ast\ast})]+(E[X_{1}])^2E[f(\mathbf{X})]\\
&+2E[X_{1}]\sum_{i=1}^{n}Cov(X_{1},X_{i})E[\nabla_{i}f(\mathbf{X}^{\ast})]+Cov(X_{1},X_{1})E[f(\mathbf{X}^{\ast})],
\end{align*}
 as required.

\section{Equivalence of the two expressions}

We shall now establish the equivalence of the two expressions in (\ref{(P11)}) and (\ref{(PP1)}) under the condition that ${\bf \Sigma}$ is positive definite. For this purpose, we will use the following lemma; see, for example, Fang et al. (1990).
 \begin{lemma} For any non-negative measurable function $f : {\Bbb R} \rightarrow {\Bbb R^+}$, we have
 $$\int_{\Bbb{R}^n} f\left(\frac{1}{2}{\bf y}^{T}{\bf y}\right)d {\bf y}=\frac{(2\pi)^{n/2}} {\Gamma(n/2)}\int_{0}^{\infty}u^{n/2-1}f(u)\mathrm{d}u.$$
 \end{lemma}

\begin{proposition}
 Under the conditions in Theorem 1, we have
 $$\frac{c_{n}}{c_{n}^{\ast}}=-\phi'(0), \;\; \frac{c_n^{\ast}}{c_{n}^{\ast\ast}}=-{\phi^{\ast}}'(0),$$
 where  $\phi$ and $\phi^{\ast}$ are   the characteristic generating functions corresponding to  the density generators  $g_n$ and $\overline{G}_{n}$, respectively.
\end{proposition}
\noindent Proof. Let ${\bf Y}=\sqrt{-\phi'(0){\bf \Sigma}^{-1}}({\bf X}-E({\bf X}))$. Then, $\mathbf{X}\sim E_{n}(\boldsymbol{\mu},~\mathbf{\Sigma},~g_{n})$
implies $\mathbf{Y}\sim E_{n}(\boldsymbol{0}, \mathbf{I_n}, g_{n})$ and $Cov(\mathbf{Y})=-\phi'(0) \mathbf{I_n}$. It then follows that
\begin{eqnarray*}
E({\bf Y}^{T}{\bf Y})&=& -n\phi'(0)\\
&=&c_n\int_{\Bbb{R}^n}{\bf y}^{T}{\bf y}g_n\left(\frac{1}{2}{\bf y}^{T}{\bf y}\right)d {\bf y} \\
&=& 2c_n \frac{(2\pi)^{n/2}} {\Gamma(n/2)}\int_{0}^{\infty}t^{n/2}g_{n}(t)\mathrm{d}t,
\end{eqnarray*}
by Lemma 1. Thus,
\begin{equation}\label{P14}
c_n=\frac{ -n\phi'(0)}{2\frac{(2\pi)^{n/2}} {\Gamma(n/2)}\int_{0}^{\infty}t^{n/2}g_{n}(t)\mathrm{d}t}.
\end{equation}
Now, by using Eq.(9) of Landsman et al. (2013), we have
$$1=c_n^*\int_{\Bbb{R}^n}  \overline{G}_{n} \left(\frac{1}{2}{\bf y}^{T}{\bf y}\right)d {\bf y}=c_n^* \frac{(2\pi)^{n/2}} {\Gamma(n/2+1)}\int_{0}^{\infty}t^{n/2}g_{n}(t)\mathrm{d}t,$$
and hence
\begin{equation}\label{P15}
c_n^*=\frac{1}{\frac{(2\pi)^{n/2}} {\Gamma(n/2+1)}\int_{0}^{\infty}t^{n/2}g_{n}(t)\mathrm{d}t}.
\end{equation}
The result that $\frac{c_{n}}{c_{n}^{\ast}}=-\phi'(0)$ readily follows from (\ref{P14}) and (\ref{P15}).

For  $\mathbf{X^*}\sim E_{n}(\boldsymbol{\mu},~\mathbf{\Sigma}, \overline{G}_{n})$, let ${\bf Z}=\sqrt{-{\phi^*}'(0){\bf \Sigma}^{-1}}({\bf X^*}-E({\bf X}^*))$. Then,
  $\mathbf{Z}\sim E_{n}(\boldsymbol{0}, \mathbf{I_n},  \overline{G}_{n})$ and $Cov(\mathbf{Z})=-{\phi^*}'(0) \mathbf{I_n}$.
By using the same arguments as above, we have
\begin{eqnarray*}
E({\bf Z}^{T}{\bf Z})&=& -n{\phi^*}'(0)\\
&=&c^*_n\int_{\Bbb{R}^n}{\bf z}^{T}{\bf z}  \overline{G}_{n}\left(\frac{1}{2}{\bf z}^{T}{\bf z}\right)d {\bf z} \\
&=& 2c^*_n \frac{(2\pi)^{n/2}} {\Gamma(n/2)}\int_{0}^{\infty}t^{n/2} \overline{G}_{n}(t)\mathrm{d}t\\
&=&nc^*_n \frac{(2\pi)^{n/2}} {\Gamma(n/2)}\int_{0}^{\infty}t^{n/2-1}  \overline{\mathcal{G}}_{n}(t)\mathrm{d}t\\
&=& \frac{nc_n^*}{c_n^{**}},
\end{eqnarray*}
and hence
$$\frac{c_n^*}{c_n^{**}}=-{\phi^*}'(0),$$
as required.

\begin{remark} From Proposition 1, we find that, when ${\bf \Sigma}$ is a positive definite, the  expressions in (\ref{(P11)}) and (\ref{(PP1)}) are indeed equivalent.
Moreover, for a positive
semidefinite matrix ${\bf \Sigma}$, we can rewrite (\ref{(PP1)}) in terms of characteristic  generators as follows:
\begin{align}\label{P16}
\nonumber E[X_{1}^{2}f(\boldsymbol{X})]=&\phi'(0){\phi^*}'(0)\sum_{i=1}^{n}\sum_{j=1}^{n}\sigma_{1i}\sigma_{1j}E[\nabla_{i,j}f(\boldsymbol{X}^{\ast\ast})]\\
&-2\mu_1\phi'(0)\sum_{i=1}^{n}\sigma_{1i}E[\nabla_{i}f(\mathbf{X}^{\ast})]\nonumber\\
&-\phi'(0)\sigma_{11}E[f({\bf X}^{\ast})]+\mu_{1}^{2}E[f(\boldsymbol{X})],
\end{align}
where  $\phi$ and $\phi^{\ast}$ are   the characteristic  generators corresponding to  the density generators  $g_n$ and $\overline{G}_{n}$, respectively.
\end{remark}

\section{Product  moments of correlated normal random variables}
 We derive now expressions for $E[X_{1}^{2}f(\boldsymbol{X})]$ of multivariate normal distribution and moments of products of correlated normal distribution.

\noindent $\mathbf{Corollary~2}$. Suppose $\mathbf{X}\sim E_{n}(\boldsymbol{\mu}, {\bf \Sigma}, \phi)$ with  density generator  $g(u)=\exp\{-u\}$ and
 characteristic generator $\phi(t)=\exp\{-t\}$. Note in this case that $g(u)=\overline{G}(u)=\overline{\mathcal{G}}(u)=\exp\{-u\}$,  $\phi(t)=\phi^*(t)$ and $b_{n}^{\ast}=b_{n}^{\ast\ast}=1$.
  Then, we have
 \begin{align*}
\nonumber E[X_{1}^{2}f(\boldsymbol{X})]=&\sum_{i=1}^{n}\sum_{j=1}^{n}\sigma_{1i}\sigma_{1j}E[\nabla_{i,j}f(\boldsymbol{X})]+2\mu_1\sum_{i=1}^{n}\sigma_{1i}E[\nabla_{i}f(\mathbf{X})]\nonumber\\
&+\sigma_{11}E[f({\bf X})]+\mu_{1}^{2}E[f(\boldsymbol{X})].\nonumber
\end{align*}
In general, for any $p_1\ge 2$, we have  the following recursive relation:
\begin{align*}
\nonumber &E[X_{1}^{p_1}f(\boldsymbol{X})]=\sum_{i=1}^{n}\sum_{j=1}^{n}\sigma_{1i}\sigma_{1j}E[\nabla_{i,j}\left(X_1^{p_1-2}f(\boldsymbol{X})\right)]\\
&+2\mu_1\sum_{i=1}^{n}\sigma_{1i}E[\nabla_{i}\left(X_1^{p_1-2}f(\mathbf{X})\right]+\sigma_{11}E[X_1^{p_1-2}f({\bf X})]+\mu_{1}^{2}E[X_1^{p_1-2}f(\boldsymbol{X})]\\
&=\sum_{i=1}^{n}\sum_{j=1}^{n}\sigma_{1i}\sigma_{1j}E[X_1^{p_1-2}\nabla_{i,j}f(\boldsymbol{X})]+\sigma_{11}^{2}(p_{1}-2)(p_{1}-3)E[X_{1}^{p_{1}-4}f(\boldsymbol{X})]\nonumber\\
&~~+2\mu_1\sum_{i=1}^{n}\sigma_{1i}E[X_1^{p_1-2}\nabla_{i}f(\mathbf{X})]+2\mu_{1}\sigma_{11}(p_{1}-2)E[X_{1}^{p_{1}-3}f(\boldsymbol{X})]\nonumber\\
&~~+\sigma_{11}E[X_1^{p_1-2}f({\bf X})]+\mu_{1}^{2}E[X_1^{p_1-2}f(\boldsymbol{X})].\nonumber
\end{align*}
This formula can be seen as a supplement to the following multivariate version of Stein's identity (see Stein (1981) and Liu (1994)):
 $$Cov(X_1,f(\boldsymbol{X}))=\sum_{i=1}^n Cov(X_1,X_i)E(\nabla_{i}f(\mathbf{X})),$$
 or equivalently,
 $$E(X_1 f(\boldsymbol{X}))=\sum_{i=1}^n \sigma_{1i}E(\nabla_{i}f(\mathbf{X}))+E(X_1)E(f(\mathbf{X})).$$
 Stein's  identity for multivariate elliptical distributions also has a similar result, that can be found in Landsman et al. (2013).

In particular, when $\phi(t)=\exp\{-t\}$,  then ${\bf{X}}\sim N_n ({\boldsymbol \mu},{\bf \Sigma})$,  ${\bf{X}^*}\sim N_n ({\boldsymbol \mu},{\bf \Sigma})$ and  ${\bf{X}^{**}}\sim N_n ({\boldsymbol \mu},{\bf \Sigma})$. In this case, we obtain  the following recursion formula:
\begin{align*}
\nonumber &E\left[\prod_{i=1}^{n}X_{i}^{p_{i}}\right]=\sigma_{11}E\left[{(X_{1})}^{p_{1}-2}\prod_{i=2}^{n}{(X_{k})}^{p_{k}}\right]+\mu_{1}^{2}E\left[X_{1}^{p_{1}-2}\prod_{k=2}^{n}X_{k}^{p_{k}}\right]\\
&~~+\sigma_{11}^{2}(p_{1}-2)(p_{1}-3)E\left[{(X_{1})}^{p_{1}-4}\prod_{k=2}^{n}{(X_{k})}^{p_{k}}\right]\\
&~~+2\sum_{j=2}^{n}\sigma_{11}\sigma_{j1}(p_{1}-2)p_{j}E\left[{(X_{1})}^{p_{1}-3}{(X_{j})}^{p_{j}-1}\prod_{k=2,~k\neq j}^{n}{(X_{k})}^{p_{k}}\right]\\
&~~+\sum_{j=2}^{n}\sigma_{j1}^{2}p_{j}(p_{j}-1)E\left[{(X_{1})}^{p_{1}-2}{(X_{j})}^{p_{j}-2}\prod_{k=2,~k\neq j}^{n}{(X_{k})}^{p_{k}}\right]\\
&~~+\sum_{j=2}^{n}\sum_{i=2,~i\neq j}^{n}\sigma_{j1} \sigma_{i1}p_{j}p_{i}E\left[{(X_{i})}^{p_{i}-1}{(X_{j})}^{p_{j}-1}{(X_{1})}^{p_{1}-2}\prod_{k=2,~k\neq i,~j}^{n}{(X_{k})}^{p_{k}}\right]\\
&~~+2\mu_{1}\bigg\{\sigma_{11}(p_{1}-2)E\left[{(X_{1})}^{p_{1}-3}\prod_{k=2}^{n}{(X_{k})}^{p_{k}}\right]\\
&~~+\sum_{j=2}^{n}\sigma_{j1}p_{j}E\left[{(X_{j})}^{p_{j}-1}{(X_{1})}^{p_{1}-2}\prod_{k=2,~k\neq j}^{n}{(X_{k})}^{p_{k}}\right]\bigg\}.
\end{align*}
\section{Results for some special cases}
In this section, we present the results for the special case of multivariate Student-$t$, logistic and Laplace distributions.\\
\noindent$\mathbf{Example~7.1}$ Multivariate Student-$t$ distribution. A $n$-dimensional Student-$t$ random vector $\mathbf{X}$, with location parameter $\boldsymbol{\mu}$, scale matrix $\mathbf{\Sigma}$ and $p>0$ degrees of freedom,  has its density function as
 \begin{align*}
f_{\boldsymbol{X}}(\boldsymbol{x})=&\frac{c_{n}}{\sqrt{|\mathbf{\Sigma}|}}\left[1+\frac{(\boldsymbol{x}-\boldsymbol{\mu})^{T}\mathbf{\Sigma}^{-1}(\boldsymbol{x}-\boldsymbol{\mu})}{p}\right]^{-\frac{p+n}{2}},~\boldsymbol{x}\in\mathbb{R}^{n},
\end{align*}
 where $c_{n}=\frac{\Gamma\left(\frac{p+n}{2}\right)}{\Gamma(p/2)(p\pi)^{\frac{n}{2}}}$. We denote it by $\mathbf{X}\sim St_{n}\left(\boldsymbol{\mu},~\boldsymbol{\Sigma},~p\right)$. In this case,
the density generator is
 $$g_{n}(t)=\left(1+\frac{2t}{p}\right)^{-(p+n)/2},$$ and so $\overline{G}_{n}(t)$ and $\overline{\mathcal{G}}_{n}(t)$ can be expressed, respectively, as
 $$\overline{G}_{n}(t)=\frac{p}{p+n-2}\left(1+\frac{2t}{p}\right)^{-(p+n-2)/2}$$
  and
 $$\overline{\mathcal{G}}_{n}(t)=\frac{p}{p+n-2}\frac{p}{p+n-4}\left(1+\frac{2t}{p}\right)^{-(p+n-4)/2}.$$ In addition, \begin{align*}
 c_{n}^{\ast}&=\frac{(p+n-2)\Gamma(n/2)}{(2\pi)^{n/2}p}\left[\int_{0}^{\infty}t^{n/2-1}\left(1+\frac{2t}{p}\right)^{-(p+n-2)/2}\mathrm{d}t\right]^{-1}\\
 &=\frac{(p+n-2)\Gamma(n/2)}{(p\pi)^{n/2}pB(\frac{n}{2},~\frac{p-2}{2})},~if~p>2
 \end{align*} and
 \begin{align*}
 c_{n}^{\ast\ast}&=\frac{(p+n-2)(p+n-4)\Gamma(n/2)}{(2\pi)^{n/2}p^{2}}\left[\int_{0}^{\infty}t^{n/2-1}\left(1+\frac{2t}{p}\right)^{-(p+n-4)/2}\mathrm{d}t\right]^{-1}\\
 &=\frac{(p+n-2)(p+n-4)\Gamma(n/2)}{(p\pi)^{n/2}p^{2}B(\frac{n}{2},~\frac{p-4}{2})},~if~p>4,
 \end{align*}
 where $\Gamma(\cdot)$ and $B(\cdot,\cdot)$ are Gamma function and Beta function, respectively.
 Then, we have
\begin{align}\label{(P17)}
\nonumber E[X_{1}^{2}f(\mathbf{X})]&= \sigma_{11}b_{n}^{\ast}E[f(\mathbf{X}^{\ast})]+
b_{n}^{\ast\ast}\sum_{i=1}^{n}\sum_{j=1}^{n}\sigma_{i1}\sigma_{j1}E[\nabla_{i,j}f(\mathbf{X}^{\ast\ast})]\\
&+2\mu_{1}b_{n}^{\ast}\sum_{i=1}^{n}\sigma_{i1}E[\nabla_{i} f(\mathbf{X}^{\ast})]+\mu_{1}^{2}E[f(\mathbf{X})],
\end{align}
where ${\bf{X}}\sim E_n ({\boldsymbol \mu},{\bf \Sigma},g_{n})$,  ${\bf{X}^*}\sim E_n ({\boldsymbol \mu},{\bf \Sigma},\overline{G}_{n})$, ${\bf{X}^{**}}\sim E_n ({\boldsymbol \mu},{\bf \Sigma},\overline{\mathcal{G}}_{n})$,
$$
b_{n}^{\ast} =\frac{p^{2}\Gamma(\frac{p+n}{2})B(\frac{n}{2},~\frac{p-2}{2})}{(p+n-2)\Gamma(\frac{p}{2})\Gamma(\frac{n}{2})}=\frac{p^{2}}{p-2},~\mathrm{if}~p>2
$$
and
$$
b_{n}^{\ast\ast} =\frac{p\Gamma(\frac{p+n}{2})B(\frac{n}{2},~\frac{p-4}{2})}{(p+n-2)(p+n-4)\Gamma(\frac{p}{2})\Gamma(\frac{n}{2})}=\frac{p}{(p-2)(p-4)},~\mathrm{if}~p>4.
$$
$\mathbf{Example~7.2}$ Multivariate logistic distribution. The density function of a $n$-dimension logistic random vector $\mathbf{X}$, with location parameter $\boldsymbol{\mu}$ and scale matrix  $\mathbf{\Sigma}$, is given by
\begin{align*}
f_{\boldsymbol{X}}(\boldsymbol{x})=\frac{c_{n}}{\sqrt{|\mathbf{\Sigma}|}}\frac{\exp\left\{-\frac{1}{2}(\boldsymbol{x}-\boldsymbol{\mu})^{T}\mathbf{\Sigma}^{-1}(\boldsymbol{x}-\boldsymbol{\mu})\right\}}{\left[1+\exp\left\{-\frac{1}{2}(\boldsymbol{x}-\boldsymbol{\mu})^{T}\mathbf{\Sigma}^{-1}(\boldsymbol{x}-\boldsymbol{\mu})\right\}\right]^{2}},~\boldsymbol{x}\in\mathbb{R}^{n},
\end{align*}
where
\begin{align*}
c_{n}&=\frac{\Gamma(n/2)}{(2\pi)^{n/2}}\left[\int_{0}^{\infty}t^{n/2-1}\frac{\exp(-t)}{[1+\exp(-t)]^{2}}\mathrm{d}t\right]^{-1}\\
&=\frac{1}{(2\pi)^{n/2}\Psi_{2}^{\ast}(-1,\frac{n}{2},1)}.
\end{align*}
Here $\Psi_{\mu}^{\ast}(z,s,a)$ is the generalized Hurwitz-Lerch zeta function defined by (cf.
Lin et al.(2006))
$$\Psi_{\mu}^{\ast}(z,s,a)=\frac{1}{\Gamma(\mu)}\sum_{n=0}^{\infty}\frac{\Gamma(\mu+n)}{n!}\frac{z^{n}}{(n+a)^{s}},$$
which has an integral representation
$$\Psi_{\mu}^{\ast}(z,s,a)=\frac{1}{\Gamma(s)}\int_{0}^{\infty}\frac{t^{s-1}e^{-at}}{(1-ze^{-t})^{\mu}}\mathrm{d}t,$$
where $\mathcal{R}(a)>0$; $\mathcal{R}(s)>0$ when $|z|\leq1~(z\neq1)$; $\mathcal{R}(s)>1$ when $z=1$.
 We denote it by $\mathbf{X}\sim Lo_{n}\left(\boldsymbol{\mu},~\boldsymbol{\Sigma}\right)$. The density generator in this case is
$$g_{n}(t)=\frac{\exp(-t)}{[1+\exp(-t)]^{2}},$$
 and  $\overline{G}_{n}(t)$ and $\overline{\mathcal{G}}_{n}(t)$ are given by
 $$\overline{G}_{n}(t)=\frac{\exp(-t)}{1+\exp(-t)},\;\;
 \overline{\mathcal{G}}_{n}(t)=\ln\left[1+\exp(-t)\right].$$ In addition,
 \begin{align*}
 c_{n}^{\ast}&=\frac{\Gamma(n/2)}{(2\pi)^{n/2}}\left[\int_{0}^{\infty}t^{n/2-1}\frac{\exp(-t)}{1+\exp(-t)}\mathrm{d}t\right]^{-1}\\
 &=\frac{1}{(2\pi)^{n/2}\Psi_{1}^{\ast}(-1,\frac{n}{2},1)}
\end{align*}
 and
 \begin{align*}
 c_{n}^{\ast\ast}&=\frac{\Gamma(n/2)}{(2\pi)^{n/2}}\left\{\int_{0}^{\infty}t^{n/2-1}\ln\left[1+\exp(-t)\right]\mathrm{d}t\right\}^{-1}\\
 &=\frac{\Gamma(n/2)}{(2\pi)^{n/2}}\left[\frac{2}{n}\int_{0}^{\infty}t^{n/2}\frac{e^{-t}}{1+e^{-t}}\mathrm{d}t\right]^{-1}\\
 &=\frac{1}{(2\pi)^{n/2}\Psi_{1}^{\ast}(-1,\frac{n}{2}+1,1)}.
\end{align*}
 Then, we have
\begin{align}\label{(P18)}
\nonumber E[X_{1}^{2}f(\mathbf{X})]&= \sigma_{11}b_{n}^{\ast}E[f(\mathbf{X}^{\ast})]+
b_{n}^{\ast\ast}\sum_{i=1}^{n}\sum_{j=1}^{n}\sigma_{i1}\sigma_{j1}E[\nabla_{i,j}f(\mathbf{X}^{\ast\ast})]\\
&+2\mu_{1}b_{n}^{\ast}\sum_{i=1}^{n}\sigma_{i1}E[\nabla_{i} f(\mathbf{X}^{\ast})]+\mu_{1}^{2}E[f(\mathbf{X})],
\end{align}
where ${\bf{X}}\sim E_n ({\boldsymbol \mu},{\bf \Sigma},g_{n})$,  ${\bf{X}^*}\sim E_n ({\boldsymbol \mu},{\bf \Sigma},\overline{G}_{n})$,  ${\bf{X}^{**}}\sim E_n ({\boldsymbol \mu},{\bf \Sigma},\overline{\mathcal{G}}_{n})$,
$$b_{n}^{\ast}=\frac{\Psi_{1}^{\ast}(-1,\frac{n}{2},1)}{\Psi_{2}^{\ast}(-1,\frac{n}{2},1)}$$
~and~
$$b_{n}^{\ast\ast}=\frac{\Psi_{1}^{\ast}(-1,\frac{n}{2}+1,1)}{\Psi_{2}^{\ast}(-1,\frac{n}{2},1)}.$$
$\mathbf{Remark~2}$ A simplification of $c_{n}$ can be found in  Yin et al. (2018):
\begin{align*}
 c_{n}=
 \begin{cases}
 &\frac{1}{\sqrt{2\pi}}[\Psi_{2}^{\ast}(-1,\frac{n}{2},1)]^{-1},~\mathrm{if}~n=1,\\
 &\frac{1}{\pi},~~~~~~~~~~~~~~~~~~~~~~~~\mathrm{if}~n=2,\\
 &\frac{1}{4\pi^{2}\ln2},~~~~~~~~~~~~~~~~~~\mathrm{if}~n=4,\\
 &\frac{1}{\pi^{n/2}(2^{n/2}-4)\zeta(\frac{n}{2}-1)},~~~~\mathrm{if}~n\geq3,~n\neq4,
 \end{cases}
 \end{align*}
 where $\zeta$ is Riemann zeta function, and it is defined as
\begin{align*}
 \zeta(s)=
 \begin{cases}
 &\sum_{n=1}^{\infty}\frac{1}{n^{s}}=\frac{1}{1-2^{-s}}\sum_{n=1}^{\infty}\frac{1}{(2n-1)^{s}},~\mathrm{if}~\mathcal{R}(s)>1,\\
 &\frac{1}{1-2^{1-s}}\sum_{n=1}^{\infty}\frac{(-1)^{n+1}}{n^{s}},~~~~\mathrm{if}~\mathcal{R}(s)>0,~s\neq1,
 \end{cases}
 \end{align*}
 which can, except for a simple pole at $s=1$ with its residue 1, be continued meromorphically to the whole complex s-plane (see Srivastava (2003) and Choi et al. (2004) for details).\\
$\mathbf{Example~7.3}$ Multivariate Laplace distribution. The density function of a Laplace random vector $\mathbf{X}$, with location parameter $\boldsymbol{\mu}$ and scale matrix  $\mathbf{\Sigma}$, is given by
 \begin{align*}
f_{\boldsymbol{X}}(\boldsymbol{x})=&\frac{c_{n}}{\sqrt{|\mathbf{\Sigma}|}}\exp\left\{-[(\boldsymbol{x}-\boldsymbol{\mu})^{T}\mathbf{\Sigma}^{-1}(\boldsymbol{x}-\boldsymbol{\mu})]^{1/2}\right\},~\boldsymbol{x}\in\mathbb{R}^{n},
\end{align*}
where $c_{n}=\frac{\Gamma(n/2)}{2\pi^{n/2}\Gamma(n)}$. We denote it by $\mathbf{X}\sim La_{n}\left(\boldsymbol{\mu},~\boldsymbol{\Sigma}\right)$. The density generator in this case is
$g_{n}(t)=\exp(-\sqrt{2t})$, and so
 $$\overline{G}_{n}(t)=(1+\sqrt{2t})\exp(-\sqrt{2t}),$$
 $$\overline{\mathcal{G}}_{n}(t)=(3+2t+3\sqrt{2t})\exp(-\sqrt{2t}).$$
 In addition,
$$c_{n}^{\ast}=\frac{n\Gamma(n/2)}{2\pi^{n/2}\Gamma(n+2)},\;\;
 c_{n}^{\ast\ast}=\frac{n(n+2)\Gamma(n/2)}{2\pi^{n/2}\Gamma(n+4)}.$$
  Then, we have
\begin{align}\label{(P19)}
\nonumber E[X_{1}^{2}f(\mathbf{X})]&= \sigma_{11}b_{n}^{\ast}E[f(\mathbf{X}^{\ast})]+
b_{n}^{\ast\ast}\sum_{i=1}^{n}\sum_{j=1}^{n}\sigma_{i1}\sigma_{j1}E[\nabla_{i,j}f(\mathbf{X}^{\ast\ast})]\\
&+2\mu_{1}b_{n}^{\ast}\sum_{i=1}^{n}\sigma_{i1}E[\nabla_{i} f(\mathbf{X}^{\ast})]+\mu_{1}^{2}E[f(\mathbf{X})],
\end{align}
where ${\bf{X}}\sim E_n ({\boldsymbol \mu},{\bf \Sigma},g_{n})$,  ${\bf{X}^*}\sim E_n ({\boldsymbol \mu},{\bf \Sigma},\overline{G}_{n})$,  ${\bf{X}^{**}}\sim E_n ({\boldsymbol \mu},{\bf \Sigma},\overline{\mathcal{G}}_{n})$,
$$b_{n}^{\ast}=n+1\;\;
{\rm ~and}\;\;~
b_{n}^{\ast\ast}=(n+3)(n+1).$$
\section{Concluding remarks}
In this work, we have derived expressions for the joint moments of elliptical distributions by two methods and have shown their equivalence when the dispersion matrix $\mathbf{\Sigma}$ is positive definite. We have used the result to derive expectations of products of correlated normal random variables and have also presented simplified  expressions for the joint moments of the special case of multivariate Student-$t$, logistic and Laplace distributions. It will be of interest to extend the results established have to truncated elliptical distributions along the lines of Shushi (2018), and to the family of skew-elliptical distributions proposed by Branco and Dey (2001). Work on these problems is currently under progress and hope to report the findings in a future paper.
\section*{Acknowledgments}
\noindent The research was supported by the National Natural Science Foundation of China (No.  11571198, 11701319).

\section*{References}
\bibliographystyle{model1-num-names}







\end{document}